\documentclass[10pt]{amsart}

\usepackage{amssymb}
\usepackage{latexsym}
\usepackage{amsmath}

\newtheorem{definition}{Definition}[section]
\newtheorem{theorem}{Theorem}[section]

\begin{document}

\title{Quantum permutation groups: a survey}
\author{Teodor Banica}
\address{T.B.: Department of Mathematics, Paul Sabatier University, 118 route de Narbonne, 31062 Toulouse,
France} \email{banica@picard.ups-tlse.fr}
\author{Julien Bichon}
\address{J.B.: Department of Mathematics, Blaise Pascal University, Campus des Cezeaux, 63177 Aubiere Cedex, France}
\email{Julien.Bichon@math.univ-bpclermont.fr}
\author{Beno\^it Collins}
\address{B.C.: Department of Mathematics,
Claude Bernard University, 43 bd du 11 novembre 1918, 69622
Villeurbanne, France, and}
\address{\vskip -6mm Department of Mathematics, University of Ottawa, 585 King Edward, Ottawa, ON K1N 6N5, Canada} \email{collins@math.univ-lyon1.fr}

\begin{abstract}
This is a presentation of recent work on quantum permutation
groups. Contains: a short introduction to operator algebras and
Hopf algebras; quantum permutation groups, and their basic
properties; diagrams, integration formulae, asymptotic laws,
matrix models; the hyperoctahedral quantum group, free wreath
products, quantum automorphism groups of finite graphs, graphs
having no quantum symmetry; complex
Hadamard matrices, cocycle twists of the symmetric group, quantum
groups acting on $4$ points; remarks and comments.
\end{abstract}

\maketitle

\section{Introduction}

The idea of noncommuting coordinates goes back to Heisenberg, who
was in turn motivated by results of Balmer and Ritz-Rydberg
regarding spectra of chemical elements. Several theories emerged
from Heisenberg's work, most complete being Connes' noncommutative
geometry, where the base space is a Riemannian manifold. See
\cite{cn}.

The specific idea of using algebras of free coordinates on
algebraic groups should be attributed to Brown \cite{br}. The
point is the following: given a group $G\subset U_n$, the matrix
coordinates $u_{ij}\in C(G)$ commute with each other, and satisfy
certain relations $R$. One can define then the universal algebra
generated by abstract variables $u_{ij}$, subject to the relations
$R$. The spectrum of this algebra is an abstract object, called
noncommutative version of $G$. The noncommutative version is not
unique, because it depends on $R$.

A detailed study of Brown's algebras, from a K-theoretic point of
view, is due to McClanahan \cite{mc}. Unfortunately, the whole
subject is a bit limited, because Brown's choice for the relations
$R$ is somehow minimal, and this makes the corresponding algebra
too big. This algebra has a comultiplication and a counit, but no
antipode. In other words, the corresponding noncommutative version
is a quantum semigroup.

The continuation of story makes use of Woronowicz's axiomatization
of compact quantum groups \cite{wo1}, \cite{wo2}. The algebras
$A_o(n)$ and $A_u(n)$, corresponding to the orthogonal and unitary
groups, appeared in Wang's thesis \cite{wa1}. Then Connes
suggested use of symmetric groups, and the algebra $A_s(n)$ was
constructed in \cite{wa2}. In all three cases the idea is the same
as Brown's. The point is to carefully choose the relations $R$, in
order to get a compact quantum group in the sense of Woronowicz.

The spectrum of $A_s(n)$ is called free version of $S_n$. This a compact quantum group, bigger than $S_n$. Its subgroups are called quantum permutation groups.

In this paper we present a number of known facts about such
quantum groups. We focus on combinatorial aspects, and their
algebraic or probabilistic interpretation.

{\bf Acknowledgements.} This work was started at Bedlewo in
October 2006 at the workshop ``Non-commutative harmonic analysis
with applications to probability''. We would like to thank Marek
Bo\.zejko for the invitation, and for several stimulating
discussions.

\section{Operator algebras}

The operator algebra background needed in order to construct
quantum permutation groups reduces to the definition of
$C^*$-algebras, and to some early work on the subject. We thought
it useful to include a short presentation of this material.
Actually the present text is written as to be at the same time an
introduction and survey paper.

\begin{definition}
A $C^*$-algebra is a complex algebra with unit, having a norm and
an involution, such that Cauchy sequences converge, and such that
$||aa^*||=||a||^2$.
\end{definition}

The basic example is $B(H)$, the algebra of bounded operators on a
Hilbert space $H$. The GNS theorem states that any $C^*$-algebra
appears as sublagebra of some $B(H)$.

The key example is $C(X)$, the algebra of continuous functions on
a compact space $X$. The Gelfand theorem below states that any
commutative $C^*$-algebra is of this form.

We need some basic spectral theory. The spectrum of an element $a\in A$ is the set $\sigma(a)$ consisting of complex numbers $\lambda$ such that $a-\lambda$ is not invertible. The spectral radius $\rho (a)$ is the radius of the smallest disk centered at $0$ containing $\sigma(a)$.

\begin{theorem}
Let $A$ be a $C^*$-algebra.
\begin{enumerate}
\item The spectrum of a norm one element is in the unit disk.
\item The spectrum of a unitary element $(a^*=a^{-1}$) is on the
unit circle. \item The spectrum of a self-adjoint element
($a=a^*$) consists of real numbers. \item The spectral radius of a
normal element ($aa^*=a^*a$) is equal to its norm.
\end{enumerate}
\end{theorem}

The first assertion follows from the formula
$1/(1-x)=1+x+x^2+\ldots$. 

If $f$ is a rational function having poles outside $\sigma(a)$, we
have $\sigma(f(a))=f(\sigma(a))$. By using the functions $z^{-1}$
and $(z+it)/(z-it)$ we get the middle assertions.

Finally, the inequality $\rho(a)\leq ||a||$ is clear from the
first assertion. For the converse we fix $\rho>\rho(a)$, and we
integrate over the circle of radius $\rho$:
$$\int
\frac{z^n}{z -a}\,dz =\sum_{k=0}^\infty\left(\int
z^{n-k-1}dz\right) a^k=a^n$$

By applying the norm and taking $n$-th
roots we get $\rho\geq\lim {||a^n||^{1/n}}$.

In the case $a=a^*$ we have $||{a^n}||=||{a}||^n$ for any exponent of the form $n=2^k$, and by taking $n$-th roots we get $\rho\geq ||{a}||$. This gives the missing inequality $\rho(a)\geq ||a||$.

In the general case $aa^*=a^*a$ we have $a^n(a^n)^*=(aa^*)^n$, and
we get $\rho(a)^2=\rho(aa^*)$. Now since $aa^*$ is self-adjoint,
we get $\rho(aa^*)=||{a}||^2$, and we are done.

\begin{theorem}
The commutative $C^*$-algebras are those of form $C(X)$.
\end{theorem}

The proof is as follows. Given a commutative $C^*$-algebra $A$, we
can define $X$ to be the set of characters $\chi :A\to\mathbb C$,
with topology making continuous all evaluation maps $e_a$. Then
$X$ is a compact space, and $a\to e_a$ is a morphism of algebras
$e:A\to C(X)$. We prove first that $e$ is involutive. We use the
following formula:
$$a=\frac{a+a^*}{2}-i\cdot\frac{i(a-a^*)}{2}$$

Thus it is enough to prove the equality $e_{a^*}=e_a^*$ for
self-adjoint elements $a$. But this is the same as proving that
$a=a^*$ implies that $e_a$ is a real function, which is in turn
true, because $e_a(\chi)=\chi(a)$ is an element of ${\sigma}(a)$,
contained in the reals.

Since $A$ is commutative, each element is normal, so $e$ is isometric: $||e_a|| =\rho(a)=||{a}||$.

It remains to prove that $e$ is surjective. But this follows from
the Stone-Weierstrass theorem, because $e(A)$ is a closed
subalgebra of $C(X)$, which separates points.

\section{Hopf algebras}

This is a short introduction to Hopf algebra philosphy.

In order to simplify presentation, we call comultiplication, counit
and antipode any morphisms of $C^*$-algebras of the following
type:
\begin{eqnarray*}
\Delta&:&A\to A\otimes A\cr \varepsilon&:&A\to {\mathbb C}\cr
S&:&A\to A^{op}
\end{eqnarray*}

The terminology comes from the fact that in the commutative case
$A=C(X)$, the morphism $\Delta$ is transpose to a binary
operation, or multiplication, $X\times X\to X$.

\begin{definition}
A finite Hopf algebra is a finite dimensional $C^*$-algebra, given
with a comultiplication, counit and antipode,
 satisfying the following conditions:
\begin{eqnarray*}
(\Delta\otimes id)\Delta&=&(id\otimes \Delta)\Delta\cr
(\varepsilon\otimes id)\Delta&=&id\cr
(id\otimes\varepsilon)\Delta&=&id\cr m(S\otimes
id)\Delta&=&\varepsilon(.)1\cr m(id\otimes
S)\Delta&=&\varepsilon(.)1
\end{eqnarray*}
\end{definition}

The group algebra $C^*(G)$ is the complex vector space spanned by $G$, with
product $g\cdot h=gh$, involution $g^*=g^{-1}$, and norm coming from the regular representation.

We say that $A$ is cocommutative if $\Sigma\Delta=\Delta$, where $\Sigma(a\otimes b)=b\otimes a$.

\begin{theorem}
Let $G$ be a finite group.
\begin{enumerate}
\item $C(G)$ is a commutative finite Hopf algebra, with
\begin{eqnarray*}
\Delta(\varphi)&=&(g,h)\to \varphi(gh)\cr
\varepsilon(\varphi)&=&\varphi(1)\cr S(\varphi)&=&g\to
\varphi(g^{-1})
\end{eqnarray*}
as structural maps. Any commutative finite Hopf algebra is of this
form. \item $C^*(G)$ is a cocommutative finite Hopf algebra, with
\begin{eqnarray*}
\Delta(g)&=&g\otimes g\cr \varepsilon(g)&=&1\cr S(g)&=&g^{-1}
\end{eqnarray*}
as structural maps. Any cocommutative finite Hopf algebra is of
this form.
\end{enumerate}
\end{theorem}

In this statement the fact that $\Delta,\varepsilon,S$ satisfy the
axioms is clear from definitions, in both situations. Observe that
the use of opposite algebras in needed for $C^*(G)$.

The assertion about commutative Hopf algebras follows from the
Gelfand theorem. For the remaining assertion, let $A$ be a finite
Hopf algebra, and consider its comultiplication, counit,
multiplication, unit and antipode. By taking duals, we get linear
maps as follows:
\begin{eqnarray*}
\Delta^*&:&A^*\otimes A^*\to A^*\cr \varepsilon^*&:&{\mathbb C}\to
A^*\cr m^*&:&A^*\to A^*\otimes A^*\cr u^*&:&A^*\to\mathbb C\cr
S^*&:&A^*\to A^*
\end{eqnarray*}

It is routine to check that these maps make $A^*$ into a finite
Hopf algebra. Moreover, if $A$ is cocommutative then $A^*$ is
commutative, so we can apply the first result. We get $A^*=C(G)$
for a certain finite group $G$, which in turn gives $A=C^*(G)$.

\section{Compact quantum groups}

There are several types of compact quantum groups. It is beyond
purposes of this paper to discuss the whole picture. Let us just
mention that we have the following list:
\begin{enumerate}
\item $q=\Omega$: subfactor quantum groups. \item $q=\sqrt{1}$:
quantum groups at roots of unity.\item $q\to 0$: asymptotic
quantum groups.\item $q>0$: deformations with positive
parameter.\item $q=1$: basic quantum groups.
\end{enumerate}

The last case is the one we are interested in. A first study here
is due to Enock and Schwartz \cite{es}. In this paper we use an
adaptation of Woronowicz's axioms in \cite{wo1}.

\begin{definition}
A finitely generated Hopf algebra is a $C^*$-algebra $A$, given
with a unitary matrix $u\in M_n(A)$ whose coefficients generate
$A$, such that the formulae
\begin{eqnarray*} \Delta(u_{ij})&=&\sum u_{ik}\otimes
u_{kj}\cr \varepsilon(u_{ij})&=&\delta_{ij}\cr
S(u_{ij})&=&u_{ji}^*
\end{eqnarray*}
define a comultiplication, a counit and an antipode.
\end{definition}

The maps $\Delta$ and $\varepsilon$ satisfy the usual axioms for a
comultiplication and a counit. The map $S$ satisfies the usual
axioms for an antipode, on the dense $*$-algebra generated by
entries of $u$. Observe also that the $q=1$ condition, namely
$S^2=id$, is satisfied.

Once the pair $(A,u)$ is given, the maps $\Delta,\varepsilon,S$
can exist or not. If they exist, they are uniquely determined, and
we have a Hopf algebra. This point of view, somehow opposite to
the spirit of abstract group theory, was invented by Woronowicz
\cite{wo1}.

The terminology and axioms are motivated by the following result.

\begin{theorem}
The following are finitely generated Hopf algebras:
\begin{enumerate}
\item $C(G)$, with $G\subset U_n$ compact Lie group. \item
$C^*(G)$, with $F_n\to G$ finitely generated group.
\end{enumerate}
\end{theorem}

In both cases, we have to exhibit a certain matrix $u$. For the
first assertion, we can use the matrix $u=(u_{ij})$ formed by
matrix coordinates of $G$, given by:
$$g=\begin{pmatrix}u_{11}(g)&&u_{1n}(g)\cr&\ddots&\cr u_{n1}(g)&&u_{nn}(g)\end{pmatrix}$$

The second assertion is clear by using the diagonal matrix
 formed by generators:
$$u=\begin{pmatrix}g_1&&0\cr&\ddots&\cr0&&g_n\end{pmatrix}$$

The algebras in the above statement can be characterized as being
the commutative or cocommutative finitely generated Hopf algebras.
See Woronowicz \cite{wo1}.

In the general case we have the following heuristic formulae:
\begin{enumerate}
\item $A=C(G)$, with $G$ compact quantum group. \item $A=C^*(G')$,
with $G'$ discrete quantum group.
\end{enumerate}

Needless to say, the quantum groups $G,G'$ don't exist as concrete
objects. This is in fact the case with all kinds of quantum
groups. See Drinfeld \cite{dr}.

\section{Free quantum groups}

We construct now the orthogonal, unitary and symmetric quantum
groups, following Wang's papers \cite{wa1}, \cite{wa2}. Let $u\in
M_n(A)$ be a square matrix over a $C^*$-algebra.
\begin{enumerate}
\item $u$ is called orthogonal if $u=\bar{u}$ and $u^t=u^{-1}$.
\item $u$ is called biunitary if $u^*=u^{-1}$ and
$u^t=\bar{u}^{-1}$.
\end{enumerate}

For the algebras $C(O_n)$ and $C(U_n)$, the corresponding matrix
$u$ is orthogonal, respectively biunitary. In the symmetric group
case the situation is less obvious. When using the embedding
$S_n\subset U_n$ given by permutation matrices, the functions
$u_{ij}$ are:
$$u_{ij}=\chi\{\sigma\in S_n\mid \sigma(j)=i\}$$

These characteristic functions satisfy a condition which reminds
magic squares:

\begin{definition}
$u\in M_n(A)$ is called magic unitary if all entries $u_{ij}$ are
projections, and on each row and column of $u$ these projections
are orthogonal, and sum up to $1$.
\end{definition}

With these definitions in hand, it is routine to check that we
have the following equalities, where $C^*_{com}$ means universal
commutative $C^*$-algebra:
\begin{eqnarray*}
C(O_n)&=&C^*_{com}\left( u_{ij}\mid u=n\times n\mbox{ {\rm
orthogonal}}\right)\cr C(U_n)&=&C^*_{com}\left( u_{ij}\mid
u=n\times n\mbox{ {\rm biunitary}}\right)\cr C(S_n)&=&C^*_{com}
\left( u_{ij}\mid u=n\times n\mbox{ {\rm magic unitary}}\right)
\end{eqnarray*}

In other words, orthogonality, biunitarity and magic unitarity are
the relevant conditions about matrix coordinates of $O_n,U_n,S_n$.
We can proceed now with liberation.

\begin{theorem}
The universal algebras
\begin{eqnarray*}
A_o(n)&=&C^*\left( u_{ij}\mid u=n\times n\mbox{ {\rm
orthogonal}}\right)\cr A_u(n)&=&C^*\left( u_{ij}\mid u=n\times
n\mbox{ {\rm biunitary}}\right)\cr A_s(n)&=&C^* \left( u_{ij}\mid
u=n\times n\mbox{ {\rm magic unitary}}\right)
\end{eqnarray*}
are finitely generated Hopf algebras.
\end{theorem}

The proof is as follows. Let us use the generic term ``special''
for the above unitarity notions. Consider now the following three
matrices, having coefficients in the target algebras of the maps
$\Delta,\varepsilon,S$ to be constructed:
\begin{eqnarray*} (\Delta u)_{ij}&=&\sum u_{ik}\otimes
u_{kj}\cr (\varepsilon u)_{ij}&=&\delta_{ij}\cr
(Su)_{ij}&=&u_{ji}^*
\end{eqnarray*}

The matrix $\varepsilon u=1$ is special, and it is routine to
check that $\Delta u$ and $Su$ are special as well. Thus the maps
$\varphi=\Delta,\varepsilon,S$ can be defined by
$\varphi(u_{ij})=(\varphi u)_{ij}$.

Summarizing, we have now free analogues of $O_n,U_n,S_n$. Their
construction might seem quite mysterious, and indeed so it is:
free quantum groups are not axiomatized.

The orthogonal and unitary algebras have the following properties:
\begin{enumerate}
\item $A_o(2)$ corresponds to the quantum group $SU_2^{-1}$. \item
$A_o(n)$ corresponds to an $R$-matrix quantization of $SU_2$.\item
$A_u(n)$ embeds into the free product $C(T)*A_o(n)$.
\end{enumerate}

We refer to \cite{bc1}, \cite{vv} for an updated discussion of
these results.

\section{Quantum permutation groups}

The algebra $A_s(n)$ is a free analogue of $C(S_n)$. We show now
that the corresponding compact quantum group consists indeed of
``quantum permutations''.

The permutations of $S_n$ act on points of $X=\{1,\ldots,n\}$. The
corresponding action map $(i,\sigma)\to\sigma(i)$ gives by
transposition a certain morphism $\alpha_{com}$, called coaction.
This coaction can be expressed in terms of the magic unitary
associated to $C(S_n)$:
$$\alpha_{com}(\delta_i)=\sum \delta_j\otimes u_{ji}$$

Now let $u$ be the magic unitary associated to $A_s(n)$, and
consider the linear map $\alpha$ given by the above formula. Then
$\alpha$ is a coaction, and we have the following result.

\begin{theorem}
We have the following commutative diagram:
$$\begin{matrix}
C(X)&\displaystyle{\mathop{\longrightarrow}^\alpha}&C(X)\otimes
A_s(n)\cr\cr \downarrow&&\downarrow\cr\cr
C(X)&\displaystyle{\mathop{\longrightarrow}^{\alpha_{com}}}&C(X)\otimes
C(S_n)
\end{matrix}$$

Moreover, $\alpha$ is the universal Hopf algebra coaction on $X$.
\end{theorem}

This result appeared in Wang's paper \cite{wa2}, in a slightly
different form. We refer to \cite{ba2} for a detailed discussion,
by using the magic unitarity condition.

At this point of writing, it is not clear whether quantum
permutations do exist. The question is whether the canonical map
$A_s(n)\to C(S_n)$ is an isomorphism or not.

\begin{theorem}
Quantum permutations exist starting from $n=4$. More precisely:
\begin{enumerate}
\item For $n=1,2,3$ we have $A_s(n)=C(S_n)$. \item For $n\geq 4$
the algebra $A_s(n)$ is not commutative, and infinite dimensional.
\end{enumerate}
\end{theorem}

The first assertion follows from the fact that for $n=1,2,3$, the
entries of a $n\times n$ magic unitary matrix have to commute with
each other. This is clear for $n=1$, and also for $n=2$, where the
magic unitary must be of the following special form:
$$u_p=\begin{pmatrix}p&1-p\cr 1-p&p\end{pmatrix}$$

At $n=3$ the proof is quite tricky. The idea is to use the Fourier
transform over $Z_3$. In terms of the vector $\xi$ formed by third
roots of unity, we can write $\alpha$ as follows:
$$\alpha(\xi^i)=\sum\xi^j\otimes v_{ji}$$

Now the magic unitarity condition on $u$ translates into a certain
condition on $v$, and the point is that with this new condition,
commutativity is clear. See \cite{ba2}.

At $n=4$ we use the following matrix, where $p,q$ are projections:
$$u_{pq}=\begin{pmatrix}
p&1-p&0&0\cr 1-p&p&0&0\cr 0&0&q&1-q\cr 0&0&1-q&q
\end{pmatrix}$$

This shows that the algebra $<p,q>$, which can be chosen to be not
commutative and infinite dimensional, is a quotient of $A_s(4)$.
This gives the last assertion.

The reasons why we have $A_s(4)\neq C(S_4)$ might remain quite
mysterious. In what follows we propose several explanations for
this fact.

\section{The Temperley-Lieb algebra}

We present here a first conceptual explanation for the main result
in previous section. The idea is that $n=4$ is the critical value
of the Jones index \cite{jo1}.

Consider the algebra $A_s(n)$. Its generators $u_{ij}$ are
coefficients of the coaction $\alpha$, so the matrix $u$ is a
corepresentation. The tensor powers of $u$ are defined as follows:
$$(u^{\otimes k})_{i_1\ldots i_k,j_1\ldots j_k}=u_{i_1j_1}\ldots
u_{i_kj_k}$$

The problem is to compute the Hom spaces between these
corepresentations.

\begin{definition}
The Temperley-Lieb algebra is given by
$$TL_n(k,l)={\rm span}\left\{ \begin{matrix}\cdot\,\cdot\,\cdot & \leftarrow &
    2k\,\,\mbox{points}\cr W & \leftarrow &
    k+l\mbox{ strings}\cr \cdot\,\cdot\,\cdot\,\cdot\,\cdot& \leftarrow &
    2l\,\,\mbox{points}\end{matrix}\right\}$$
    where strings join pairs of points, do not cross, and are taken up to isotopy.
    \end{definition}

$TL_n$ is a tensor category: the composition is by vertical
concatenation, with the rule that closed circles are deleted and
replaced by the number $n$, the tensor product is by horizontal
concatenation, and the involution is by upside-down turning of
diagrams.

Temperley-Lieb diagrams act on tensors according to the following
formula, where the middle symbol is $1$ if all strings of $p$ join
pairs of equal indices, and is $0$ if not:
$$p(e_{i_1}\otimes\ldots\otimes e_{i_k})=\sum_{j_1\ldots j_l}\begin{pmatrix}i_1\,i_1\ldots i_k\,i_k\cr p\cr j_1\,j_1\ldots j_l\,j_l\end{pmatrix}e_{j_1}\otimes\ldots\otimes e_{j_l}$$

In case the index satisfies $n\geq 4$, different diagrams produce
different linear maps, and this action makes $TL_n$ a subcategory
of the category of Hilbert spaces.

\begin{theorem}
We have the equality of tensor categories
$$Hom(u^{\otimes k},u^{\otimes l})=TL_n(k,l)$$
where $u$ is the fundamental corepresentation of $A_s(n)$, with
$n\geq 4$.
\end{theorem}

The proof uses Woronowicz's duality in \cite{wo2}. The idea is
that the definition of $A_s(n)$ translates into a presentation
result for the corresponding tensor category:
\begin{enumerate}\item The fact that a unitary matrix $u$ is magic is equivalent to $M\in Hom(u^{\otimes 2},u)$ and $U\in
Hom(1,u)$, where $M,U$ are the multiplication and unit of $C^n$.
\item The relations satisfied by $M,U$ in a categorical sense are
those satisfied by the diagrams $m=|\cup|$ and $u=\cap$, which in
turn generate $TL_n$.
\end{enumerate}

We recall now that the tensor powers of the fundamental
representation of $PU_2\simeq SO_3$ form a category which is
isomorphic to $TL_4$. This shows that the irreducible
representations of $A_s(4)$, with fusion rules and dimensions, are
the same as irreducible representations of $SO_3$. Moreover, the
fusion rule statement must hold for any $n\geq 4$.

\begin{theorem}
For any $n\geq 4$, the irreducible corepresentations of $A_s(n)$
satisfy the Clebsch-Gordan rules for irreducible representations
of $SO_3$.
\end{theorem}

In other words, the irreducible corepresentations are as follows:
\begin{enumerate}\item They are given by $r_0=1$, $r_1=u-1$, $r_2=u^{\otimes 2}-3u+1$ and so on.
\item They satisfy $r_a\otimes r_b=r_{|a-b|}+r_{|a-b|+1}+\ldots
+r_{a+b-1}+r_{a+b}$.
\end{enumerate}

These considerations have several extensions. We would like to
mention here the following statement, which trivializes the whole
thing: the quotients of $A_s(n)$ are in functorial correspondence
with subalgebras of the $n$-th spin planar algebra. See
\cite{ba3}.

\section{The Weingarten formula}

By general results of Woronowicz in \cite{wo1}, the Hopf algebra $A_s(n)$ has a unique unital bi-invariant state, called Haar integration, and denoted here as an integral:
$$\int :A_s(n)\to\mathbb C$$

The various integrals can be computed by using the representation theory diagrams found in previous section. The idea here, going back to Weingarten's paper \cite{we}, was developed in \cite{bia}, \cite{co}, \cite{cs1}, \cite{cs2}, and was applied to quantum groups in \cite{bc1}, \cite{bc2}, \cite{bc3}.

\begin{definition}
Consider the set $NC(k)$ of non-crossing partitions of $\{1,\ldots ,k\}$.
\begin{enumerate}
\item We plug multi-indices $i=(i_1,\ldots ,i_k)$ into partitions $p\in NC(k)$, and we set $\delta_{pi}=1$ if all blocks of $p$ contain equal indices of $i$, and $\delta_{pi}=0$ if not.
\item The Gram matrix of partitions (of index $n\geq 4$) is given by $G_{kn}(p,q)=n^{|p\vee q|}$, where $\vee$ is the set-theoretic sup, and $|.|$ is the number of blocks.
\item The Weingarten matrix $W_{kn}$ is the inverse of $G_{kn}$.
\end{enumerate}
\end{definition}

The non-crossing partitions are in correspondence with Temperley-Lieb diagrams having no upper points: these can be indeed obtained by fattening the partitions.

Now by using the Temperley-Lieb action described in previous section, we see that the elements of $NC(k)$ create a basis of fixed vectors of $u^{\otimes k}$. The Gram matrix of this basis is nothing but $G_{kn}$, as shown by the following computation:
$$<p,q>=\sum_i\delta_{pi}\delta_{qi}=\sum_i\delta_{p\vee q,i}=n^{|p\vee q|}$$

Observe also that $G_{kn}$ is by definition a kind of version of Di Francesco's meander matrix in \cite{df}. With these notations, we have the following result.

\begin{theorem}
The Haar functional of $A_s(n)$ is given by
$$\int u_{i_1j_1}\ldots u_{i_kj_k}=\sum_{pq}\delta_{pi}\delta_{qj}W_{kn}(p,q)$$
where the sum is over all pairs of diagrams $p,q\in NC(k)$.
\end{theorem}

The proof is based on the following fact: the numbers on the left are the matrix coefficients of the orthogonal projection onto the space of fixed points of $u^{\otimes k}$.

As a first consequence, we have the following moment formula:
$$\int (u_{11}+\ldots+u_{ss})^k=Tr(G_{kn}^{-1}G_{ks})$$

The free Poisson law of parameter $t\in (0,1]$ is the following 
probability measure:
$$\pi_t=(1-t)\,\delta_0+\frac{1}{2\pi x}\sqrt{4t -(x-1-t)^2}\,dx$$

This measure is also called Marchenko-Pastur law. The terminology here comes from the fact that $\pi_t$ is the free analogue of the Poisson law of parameter $t$. See \cite{vo}, \cite{vdn}.

\begin{theorem}
The variable $u_{11}+\ldots +u_{ss}$ with $s=[tn]$ and $n\to\infty$ has law $\pi_t$.
\end{theorem}

This follows from the moment formula, by using the fact that with $n\to\infty$, both the Gram and Weingarten matrices are concentrated on the diagonal. The trace to be computed reduces to a sum of powers of $t$, known to give the $k$-th moment of $\pi_t$.

In the classical case a similar result is available, in terms of Poisson laws. As a conclusion, $C(S_n)\to A_s(n)$ transforms asymptotic independence into freeness. See \cite{bc2}.

\section{The Pauli quantum group}

The central object of the theory is the algebra $A_s(4)$. In this
section we present an explicit matrix model for this algebra,
coming from the Pauli matrices:
$$c_1=\begin{pmatrix}1&0\cr 0&1\end{pmatrix}\hskip 5mm
c_2=\begin{pmatrix}i&0\cr 0&-i\end{pmatrix}\hskip 5mm
c_3=\begin{pmatrix}0&1\cr -1&0\end{pmatrix}\hskip 5mm
c_4=\begin{pmatrix}0&i\cr i&0\end{pmatrix}$$

These matrices multiply according to the formulae for quaternions:
$$c_2^2=c_3^2=c_4^2=-1$$
$$c_2c_3=-c_3c_2=c_4$$
$$c_3c_4=-c_4c_3=c_2$$
$$c_4c_2=-c_2c_4=c_3$$

The Pauli matrices form an orthonormal basis of $M_2(\mathbb C)$,
and the same is true if we multiply them to the left or to the
right by an element of $SU_2$. This shows that for any $x\in
SU_2$, the elements $\xi_{ij}=c_ixc_j$ form a magic basis of
$M_2({\mathbb C})$, in the sense that the corresponding orthogonal
projections $P_{ij}$ form a magic unitary over $M_4({\mathbb C})$.

\begin{definition}
The Pauli representation is the map
$$\pi:A_s(4)\to C(SU_2, M_4({\mathbb C}))$$
given by $\pi(u_{ij})=(x\to$ rank one projection on $c_ixc_j)$.
\end{definition}

This representation is introduced in \cite{bm}. In \cite{bc3} we
use integration techniques for proving that $\pi$ is faithful. The
idea is to check commutativity of the following diagram:
$$\begin{matrix}
A_s(4)&\ &\to&\ &C(SU_2,M_4(\mathbb C))\cr \ \cr \downarrow&\ &\
&\ &\downarrow\cr \ \cr\mathbb C&\ & \leftarrow&\ &M_4({\mathbb
C})
\end{matrix}$$

A key problem is to work out the integral geometric analogy
between $C(S_4)$ and $A_s(4)$, at level of laws of averages of
diagonal coordinates $u_{ii}$. For $C(S_4)$ we have:
$${\rm law}(t_1u_{11}+\ldots+t_4u_{44})=\frac{1}{24}\left(9\delta_0+\delta_1+2\sum_i
\delta_{t_i}+\sum_{i\neq j} \delta_{t_i+t_j}\right)$$

For the algebra $A_s(4)$ we can use the Pauli representation,
which makes integration problems correspond to computations on the
real sphere $S^3$.

\begin{theorem}
For $s=1,2,4$ we have the formula
$${\rm law}\left(s^{-1}(u_{11}+\ldots+u_{ss})\right)=\left(1-\frac{s}{4}\right)\,\delta_0+
\frac{s}{4}\,\mu_s$$ where $\mu_1,\mu_2,\mu_4$ are a Dirac mass, a
Lebesgue measure, and a free Poisson law.
\end{theorem}

This result doesn't quite clarify the relation with $C(S_4)$. In
fact, what is missing is the $s=3$ law. This is the law of the
following matrix, depending on $(a,b,c,d)\in S^3$:
$$M_3=\frac{1}{3}\begin{pmatrix}
3a^2&-ab&-ac&-ad\cr -ab&3b^2&-bc&-bd\cr -ac&-bc&3c^2&-cd\cr
-ad&-bd&-cd&3d^2
\end{pmatrix}$$

We don't know how to compute this measure. The problem is
explained in \cite{bc3}.

\section{The hyperoctahedral quantum group}

In this section we present a few known facts, along with some recent work \cite{ahn}. The quantum symmetry algebra of a finite graph $X$
is defined as follows:
\begin{enumerate}
\item In case $X$ has no edges, we set $A(X)=A_s(n)$. The fact
that this is indeed a quantum symmetry algebra follows from
considerations in previous sections. \item In the general case, we
set $A(X)=A_s(n)/R$, where $R$ are the relations coming from
$du=ud$, where $d\in M_n(0,1)$ is the adjacency matrix of $X$.
\end{enumerate}

The quotients of $A(X)$ are called Hopf algebras coacting on $X$.
See \cite{ba2}, \cite{bb1}, \cite{bi2}.

The simplest example is with the square. Let us go back to the
magic unitary matrix $u_{pq}$ in section 6. This matrix commutes
with the adjacency matrix of the square. Moreover, by choosing the
projections $p,q$ to be free, the algebra $<p,q>$ they generate is
isomorphic to the group algebra of $Z_2*Z_2=D_\infty$, and we get
the following result.

\begin{theorem}
The dual of $D_\infty$ acts on the square.
\end{theorem}

In other words, for the square we have an arrow $A(X)\to
C^*(D_\infty)$. The algebra $A(X)$ can be actually computed
explicitly, and is isomorphic to $C(O_2^{-1})$. See \cite{bi3}.

This result has the following generalization. Consider the Hopf
algebra $C(O_n^{-1})$, which is the quotient of $A_o(n)$ by the
skew-commutation relations for $GL_n^{-1}$:
\begin{enumerate}
\item  $u_{ij}u_{ik}=-u_{ik}u_{ij}$, $u_{ji}u_{ki} =
-u_{ki}u_{ji}$,  for $i\neq j$. \item $u_{ij}u_{kl}=u_{kl}u_{ij}$
for $i\neq k$, $j\neq l$.
\end{enumerate}

These relations define a Hopf ideal, so we have indeed a Hopf algebra.

\begin{theorem}
$C(O_n^{-1})$ is the quantum symmetry algebra of the hypercube in
$\mathbb R^n$.
\end{theorem}

This leads to the quite surprising conclusion that $O_n^{-1}$ is a
quantum analogue of the hyperoctahedral group $H_n$. On the other
hand, $O_n^{-1}$ cannot be a free version of $H_n$, say because
the fusion semiring depends on $n$, which avoids probabilistic
freeness.

In order to solve this problem, the idea is as follows. The group
$H_n$ appears also as symmetry group of the space formed by the
$[-1,1]$ segments on each coordinate axis. In other words, $H_n$
is the symmetry group of $I_n$, the graph formed by $n$ segments.

\begin{definition}
$A_h(n)$ is the quantum symmetry algebra of $I_n$.
\end{definition}

This algebra is the quotient of $A_s(2n)$ by the relations coming
from commutation of $u$ with the adjacency matrix of $I_n$. Now
writing down the commutation relations leads to the conclusion
that $u$ must be a magic unitary of the following special form:
$$u=\begin{pmatrix}a&b\cr b&a\end{pmatrix}$$

We call such a matrix a $2n\times 2n$ sudoku unitary. With this
notion in hand, $A_h(n)$ is the universal algebra generated by
entries of a $2n\times 2n$ sudoku unitary.

As a first consequence, we have the following commutative diagram:
$$\begin{matrix}
A_u(n)&\to&A_o(n)&\to&A_h(n)&\to&A_s(n)\cr\cr
\downarrow&&\downarrow&&\downarrow&&\downarrow\cr\cr
C(U_n)&\to&C(O_n)&\to&C(H_n)&\to&C(S_n)
\end{matrix}$$

A Tannakian translation gives the representation theory of
$A_h(n)$: the relevant algebra is the Fuss-Catalan algebra of
Bisch and Jones \cite{bj1}, and the Weingarten machinery leads to
a free analogue of the Bessel function combinatorics for $H_n$.

\section{Free wreath products}

The simplest example of a wreath product is the hyperoctahedral
group $H_n$. Consider the graph $I_n$ formed by $n$ segments, and
for each segment $I^i$ consider the element $\tau_i\in H_n$ which
returns $I^i$, and keeps the other segments fixed. The elements
$\tau_1\ldots\tau_n$ have order two and commute with each other,
so they generate a product $L_n=Z_2\times\ldots\times Z_2$. Now
the symmetric group $S_n$ acts as well on $I_n$ by permuting the
segments, and it is routine to check that we have $H_n=L_n\rtimes
S_n$. But this latter group is a wreath product:

\begin{theorem}
We have $H_n=Z_2\wr S_n$.
\end{theorem}

At level of algebras of functions, this gives a formula of the
following type, where $\times_w$ is some functional analytic
implementation of the $\wr$ operation:
$$C(H_n)=C(Z_2)\times_w C(S_n)$$

The operation $\times_w$ is not the good one for quantum groups.
This is because the natural quantum formula involving it would
imply $A_h(2)=C(H_2)$, which is wrong. The point is that in the
quantum world, wreath products are replaced by free wreath
products.

\begin{definition}
The free wreath product of $(A,u)$ and $(B,v)$ is given by
$$A*_wB=(A^{*n}*B)/<[u_{ij}^{(a)},v_{ab}]=0>$$
where $n$ is the size of $v$, and has magic unitary matrix
$w_{ia,jb}=u_{ij}^{(a)}v_{ab}$.
\end{definition}

This notion, introduced in \cite{bi3}, is justified by formulae of
the following type, where $G,A$ denote classical symmetry groups,
respectively quantum symmetry algebras:
\begin{eqnarray*}
G(X*Y)&=&G(X)\ \wr\  G(Y)\cr A(X*Y)&=&A(X)*_wA(Y)\end{eqnarray*}

There are several such formulae, and the one we are interested in
is:
\begin{eqnarray*}
G(I\ldots I)&=&G(I)\ \wr \ G(\circ\ldots\circ)\cr A(I\ldots
I)&=&A(I)*_wA(\circ\ldots\circ)\end{eqnarray*}

Here $I$ is a segment, $\circ$ is a point, and the dots mean
$n$-fold disjoint union. The upper formula is $H_n=Z_2\wr S_n$. As
for the lower formula, this is what we need:

\begin{theorem}
We have $A_h(n)=C(Z_2)*_wA_s(n)$.
\end{theorem}

Now getting back to general free wreath products, a first thing to
be noticed is the following diagram, where maps on the left are
defined by formulae on the right:
$$\begin{matrix}
A^{*n}*B&\to&A*_wB\cr \ \cr \downarrow&\ &\downarrow\cr \ \cr
A*B&\to&A\otimes B
\end{matrix}\hskip7mm :\hskip7mm\begin{matrix}
\sum u_{ii}^{(a)}v_{aa}&\to&\sum u_{ii}^{(a)}v_{aa}\cr \ \cr
\downarrow&\ &\downarrow\cr \ \cr \sum u_{ii}v_{aa}&\to&\sum
u_{ii}\otimes v_{aa}
\end{matrix}$$

The spectral measures of the north-east and south-west elements
can be computed in several cases of interest, and turn out to be
equal. The conjecture is that equality holds, under general
assumptions. This is the same as saying that the spectral measure
of $A*_wB$ appears as free multiplicative convolution of the
spectral measures of $A,B$:
$$\mu (A*_wB)=\mu (A)\boxtimes \mu(B)$$

This is related to work in preparation of Bisch and Jones, see
\cite{bb1}. The whole thing would be a first step towards
establishing an analogy with results of \'Sniady \cite{sn}.

\section{Quantum automorphisms of finite graphs}

The free wreath product results in previous section are part of a
classification project for finite graphs. This is in turn part of the
Bisch-Jones classification project for planar algebras generated
by a $2$-box \cite{bj2}, \cite{bj3}. Indeed, it follows from the
general Tannaka-Galois duality in \cite{ba3} that the algebras of
form $A(X)$ with $X$ colored oriented graph are those
corresponding to subalgebras of the spin planar algebra, generated
by a $2$-box.

In the Bisch-Jones project the complexity is measured by the
dimension of the second commutant, which is $d=15,16$. This choice
is perfectly legitimate in the context of general planar algebras.
In the our particular case however, a more natural invariant is
the number $n$ of vertices. This number is $n=7,8,9,11$ in the
papers \cite{ba2}, \cite{ba3}, \cite{bb1}, \cite{bb2}.

The conclusion in \cite{bb2} is the following table, containing
all vertex-transitive graphs of order $n\leq 11$ modulo
complementation, except for the Petersen graph:

\begin{center}\begin{tabular}[t]{|l|l|l|l|}\hline
Order&Graph&Classical group&Quantum group\\ \hline\hline 2&$K_2$ {\rm (simplex)}&$Z_2$&$Z_2$\\
\hline\hline 3&$K_3$&$S_3$&$S_3$\\ \hline\hline
4&$2K_2$ {\rm (duplication)}&$H_2$&$H_2^+$ {\rm (hyper. quant. group)}\\
\hline 4&$K_4$&$S_4$&$S_4^+$ {\rm (symm. quant. group)}\\
\hline\hline 5&$C_5$ {\rm (cycle)}&$D_5$&$D_5$\\ \hline
5&$K_5$&$S_5$&$S_5^+$\\ \hline\hline 6&$C_6$&$D_6$&$D_6$\\ \hline
6&$2K_3$&$S_3\wr Z_2$&$S_3{\,\wr_*\,}Z_2$ {\rm (free wreath prod.})\\
\hline 6&$3K_2$&$H_3$&$H_3^+$\\ \hline 6&$K_6$&$S_6$&$S_6^+$\\
\hline\hline 7&$C_7$&$D_7$&$D_7$ \\ \hline
7&$K_7$&$S_7$&$S_7^+$\\ \hline\hline 8&$C_8$, $C_8^+$ {\rm (cycle with diags.)}&$D_8$&$D_8$\\
\hline 8&$P(C_4)$ {\rm (prism)}& $H_3$&$S_4^+\times Z_2$ \\ \hline 8&$2K_4$&$S_4\wr Z_2$&$S_4^+{\,\wr_*\,}Z_2$ \\
\hline 8&$2C_4$& $H_2\wr Z_2$ & $H_2^+{\,\wr_*\,}Z_2$
\\ \hline 8&$4K_2$&$H_4$&$H_4^+$ \\ \hline 8&$K_8$&$S_8$&$S_8^+$\\
\hline\hline 9&$C_9$, $C_9^3$ {\rm (cycle with
chords)}&$D_9$&$D_9$\\ \hline 9 & $K_3\times K_3$ {\rm (discrete
torus)}&$S_3\wr Z_2$&$S_3\wr Z_2$
\\ \hline 9&$3K_3$&$S_3\wr S_3$&$S_3{\,\wr_*\,}S_3$ \\ \hline
9&$K_9$&$S_9$&$S_9^+$ \\ \hline \hline 10&$C_{10}$, $C_{10}^2$,
$C_{10}^+$, $P(C_5)$&$D_{10}$&$D_{10}$\\ \hline 10 &
$P(K_5)$&$S_5\times Z_2$&$S_5^+\times Z_2$\\ \hline
10&$C_{10}^4$&$Z_2\wr D_5$&$Z_2{\,\wr_*\,}D_5$\\ \hline
10&$2C_5$&$D_5\wr Z_2$&$D_5{\,\wr_*\,}Z_2$\\ \hline
10&$2K_{5}$&$S_5\wr Z_2$&$S_5^+{\,\wr_*\,}Z_2$\\ \hline
10&$5K_2$&$H_5$&$H_5^+$\\ \hline 10&$K_{10}$&$S_{10}$&$S_{10}^+$\\
\hline\hline 11&$C_{11}$, $C_{11}^2$,
$C_{11}^3$&$D_{11}$&$D_{11}$\\ \hline
11&$K_{11}$&$S_{11}$&$S_{11}^+$\\ \hline
\end{tabular}\end{center}

\section{Graphs having no quantum symmetry}

The classification project for finite graphs is there for various
reasons, one of them being to help in classification of certain
subfactors and planar algebras, of integer index. In other words,
the whole thing should be regarded as belonging to a specialized
area of von Neumann algebras, and mathematical physics in general.

As explained to us by Jones, the end of the game would be to
investigate some special graphs, such as the Clebsch graph, or the
Higman-Sims graph.

We are quite far away from this kind of application. The big list
in previous section consists of simplexes, cycles, and various
products between them. That is, our main realization so far is to
have reasonably strong results about product operations.

The next step would be to develop some new techniques, for graphs
which do not decompose as products. The first such graph is the
Petersen one, at $n=10$. As already mentioned, we have no results
about it. But work here is in progress, and we hope to come up
soon with an answer, along with a study for higher $n$, say
between 12-15.

Summarizing, we think that the pedestrian approach in \cite{ba2},
\cite{ba3}, \cite{bb1}, \cite{bb2} must be continued with a few
more papers, in order for the theory to be mature enough.

We would like to present now a first conceptual result emerging
from our small $n$ study. This concerns graphs having no quantum
symmetry.

\begin{definition}
A finite graph $X$ has no quantum symmetry if it satisfies one of
the following equivalent conditions, where $d$ is its adjacency
matrix:
\begin{enumerate}
\item The quantum symmetry algebra $A(X)$ is commutative. \item We
have $A(X)=C(G_X)$, where $G_X$ is the symmetry group of $X$.
\item For a magic unitary $u$, $du=ud$ implies that $u_{ij}$
commute with each other.
\end{enumerate}
\end{definition}

The problem of characterizing such graphs goes back to Wang's
paper \cite{wa2}, with the results $A(C_3)=C(S_3)$ and $A(C_4)\neq
C(H_2)$, showing that $C_3$ has no quantum symmetry, but $C_4$
does. A topological formulation of the problem is found by Curtin
in \cite{cu}.

There are several graphs in the above table which satisfy this
condition: the cycles $C_n$ with $n\neq 4$, a number of prisms and
of cycles with chords, and the discrete torus $K_3\times K_3$.
Moreover, we have found some more graphs by working on the
subject, and this led us to take a detailed look at the case of
circulant graphs.

A graph $X$ having $n$ vertices is called circulant if its
automorphism group contains copy of $Z_n$. This is the same as
saying that vertices of $X$ are elements of $Z_n$, and that $i\sim
j$ (connection by an edge) implies $i+k\sim j+k$ for any $k$.

Associated to a circulant graph are the following algebraic
invariants:
\begin{enumerate}
\item The set $S\subset Z_n$ is given by $i\sim j \iff j-i \in S$.
\item The group $E\subset Z_n^*$ consists of elements $a$ such
that $aS=S$. \item The order of $E$ is denoted $k$, and is called
type of $X$.
\end{enumerate}

With these notations, we have the following result:

\begin{theorem}
A type $k$ circulant graph having $p>>k$ vertices, with $p$ prime,
has no quantum symmetry.
\end{theorem}

This result is proved in \cite{bbc}, with the lower bound
$p>6^{{\varphi(k)}}$, where $\varphi$ is the Euler function. Most
of the proof there doesn't really use the fact that $p$ is prime,
and we hope to come up soon with more general results in this
sense.

The whole subject is somehow opposite to the freeness
considerations in previous sections: ``no quantum symmetry'' means
``too many constrains, hence independence''.

\section{Hadamard matrices}

A complex Hadamard matrix is a matrix $h\in M_n(\mathbb C)$ having
the following properties: entries are on the unit circle, and rows
are mutually orthogonal.

The rows of $h$, denoted $h_1,\ldots ,h_n$, can be regarded as
elements of the algebra ${\mathbb C}^n$. Since each $h_i$ is
formed by complex numbers of modulus $1$, this element is
invertible. We can therefore consider the elements
$\xi_{ij}=h_j/h_i$, and we have:
\begin{eqnarray*}
<\xi_{ij},\xi_{ik}>&=&<h_j/h_i,h_k/h_i>\cr &=&n<h_j,h_k>\cr
&=&n\cdot\delta_{jk}
\end{eqnarray*}

A similar computation works for columns, so $\xi$ is a magic basis
of $\mathbb C^n$, in the sense that the corresponding orthogonal
projections $P(\xi_{ij})$ form a magic unitary matrix.

\begin{definition}
Let $h\in M_n(\mathbb C)$ be an Hadamard matrix. \begin{enumerate}
\item $\xi$ is the magic basis of $\mathbb C^n$ given by
$\xi_{ij}=h_j/h_i$. \item $P$ is the magic unitary over
$M_n(\mathbb C)$ given by $P_{ij}=P(\xi_{ij})$. \item $\pi$ is the
representation of $A_s(n)$ given by $\pi(u_{ij})=P_{ij}$. \item
$A$ is the quantum permutation algebra associated to $\pi$.
\end{enumerate}
\end{definition}

In other words, we say that the representation $\pi:A_s(n)\to
M_n(\mathbb C)$ comes from a representation $\nu:G_n\to U_n$ of
the dual of the $n$-th quantum permutation group, then we consider
the quantum group $G=\nu(G_n)$, and the algebra $A=C^*(G)$. See
\cite{bn}.

\begin{theorem}
The construction $h\to A$ has the following properties:
\begin{enumerate}
\item The Fourier matrix $h_{ij}=w^{ij}$ with $w=e^{2\pi i/n}$
gives $A=C(Z_n)$.\item For a tensor product $h=h'\otimes h''$ we
have $A=A'\otimes A''$.\item $A$ is commutative if and only if $h$
is a tensor product of Fourier matrices.
\end{enumerate}
\end{theorem}

For $n=1,2,3,5$ any Hadamard matrix is equivalent to the Fourier
one. At $n=4$ we have the following example, depending on $q$ on
the unit circle:
$$h_q=\begin{pmatrix}
1&1&1&1\cr 1&q&-1&-q\cr 1&-1&1&-1\cr 1&-q&-1&q
\end{pmatrix}$$

These are, up to equivalence, all $4\times 4$ Hadamard matrices.
As an example, the Fourier matrix corresponds to the value $q=\pm
i$. See Haagerup \cite{ha}.

\begin{theorem}
Let $n$ be the order of $q^2$, and for $n<\infty$ write $n=2^sm$,
with $m$ odd. The matrix $h_q$ produces the algebra $C^*(G)$,
where $G$ is as follows:
\begin{enumerate}
\item For $s=0$ we have $G=Z_{2n}\rtimes Z_2$. \item For $s=1$ we
have $G=Z_{n/2}\rtimes Z_4$.\item For $s\geq 2$ we have
$G=Z_{n}\rtimes Z_4$.\item For $n=\infty$ we have $G=Z\rtimes
Z_2$.
\end{enumerate}
\end{theorem}

This result provides the first example of a deformation situation
for quantum permutation groups. Observe that the parameter space
is the unit circle, with roots of unity highlighted. We should
mention that this space, while being fundamental in most theories
emerging from Drinfeld's original work \cite{dr}, is quite new in
the compact quantum group area, where the deformation parameter is
traditionally real. See \cite{bn}.

There are many difficult problems regarding Hadamard matrices, and
we don't know yet if quantum permutation groups can help. See
Jones \cite{jo2}.

\section{Cocycle twists of the symmetric group}

The examples of quantum permutation groups discussed so far in this
paper are either classical or infinite-dimensional. A natural question
is whether there exist non-classical finite quantum permutation groups.
The construction of such objects was done in \cite{bi1}, using the 2-cocycle
twisting procedure.

The idea of twisting, originally due to Drinfeld \cite{dr2}, and developed
by Doi \cite{do} in the dual framework, is the following one.
Starting with a Hopf algebra $A$, we consider linear maps
$\sigma : A \otimes A \to \mathbb C$ satisfying certain conditions
and called (Hopf) 2-cocycles. We  then deform the product of
$A$ by $\sigma$ to get a new Hopf algebra $A^{\sigma}$, called a twist of $A$, and having
the same tensor category of corepresentations as $A$. See Schauenburg \cite{sc}.

The theory of twisting is developed at different levels of
generality and studied in numerous papers that we shall not list
here. Amongst these, the paper of Enock and Vainerman \cite{ev}
was influential: they realized that 2-cocycles could be easily
constructed from abelian subgroups, over which Hopf 2-cocycles
correspond to ordinary group 2-cocycles. This idea was used in
\cite{bi1} to construct twists of $C(S_{2n})$ induced by the
abelian subgroup $Z_2^n$, leading to the following construction.

Let $i \in \{1,\ldots, 2n\}$. For $i$ even we put $i'= i-1$ and
$i^*= i/2$. For  $i$ odd we put $i' = i + 1$ and $i^* = i'/2$.
We consider a matrix $p=(p_{ij}) \in
  M_n(\mathbb C)$ with $p_{ii} = 1$ and $p_{ij} = p_{ji} = \pm 1$ for all
    $i$ and $j$.

\begin{definition}
The Hopf algebra $C_p(S_{2n})$ is the quotient of
$A_s(2n)$ by the following relations:
    \begin{eqnarray*}
&(3 + p_{i^*j^*})u_{kj}u_{li} + (1 - p_{i^*j^*})u_{kj}u_{li'}
+ (1 - p_{i^*j^*})u_{kj'}u_{li} +(p_{i^*j^*} -1)u_{kj'}u_{li'} \\
= &\ (3 + p_{l^*k^*})u_{li}u_{kj} + (1 - p_{l^*k^*})u_{l'i}u_{kj}
+ (1 - p_{l^*k^*})u_{li}u_{k'j} + (p_{l^*k^*} - 1)u_{l'i}u_{k'j} \
\end{eqnarray*}
\end{definition}

The conceptual meaning of these relations is that they are
FRT-relations \cite{rtf} associated with a Yang-Baxter operator
$\mathbb C^{2n} \otimes \mathbb C^{2n} \to \mathbb C^{2n} \otimes
\mathbb C^{2n}$ attached to $p$.

One shows that $C_p(S_{2n})$ is a twist of $C(S_{2n})$ and here is the conclusion in \cite{bi1}.

\begin{theorem}
There exist at least $n$ non-isomorphic
 finite quantum permutation groups acting on $2n$ points, and
 having the same tensor category of representations as $S_{2n}$.
\end{theorem}

It is also possible to construct twists of $C(S_{2n+1})$ in the same manner.
More general twistings of $S_n$ are constructed in \cite{bi20}, using
arbitrary roots of unity.

We have the following natural questions.

\begin{enumerate}
\item Does there exist a finite graph having a non-classical
finite quantum symmetry group, for example one of the above ones
constructed by twisting? \item Consider a  finite quantum group
$G$ obtained as a twisting of $S_n$. Is it true that $G$ is a
quantum permutation group?
\end{enumerate}

The twisting construction is also available for any finite group.
This leads to some surprises: although the alternating group $A_5$
does not act faithfully on $4$ points, it has a finite quantum
analogue that does. This plays an important role in the
classification of the quantum groups acting on 4 points in the
next section.

Finally we should mention that the twisting procedure is also a
useful tool to understand the infinite-dimensional situation: for
example the Hopf algebras $C(O_n^{-1})$ are twists of $C(O_n)$. As
explained in the next section, in the case $n=4$, twisting
techniques essentially enable us to classify the quantum groups
acting on $4$ points.

\section{Quantum groups acting on $4$ points}

A natural problem in the area of quantum permutation groups is the
classification problem, at least for small $n$. This is the same as the classification of Hopf algebra quotients of $A_s(n)$.
In the dual language, we have to classify the quantum
subgroups of $S_n^+$, the compact quantum group dual to $A_s(n)$,
defined by $A_s(n)=C(S_n^+)$.

At $n=4$ we have the following result \cite{bb3}.

\begin{theorem}
 The compact quantum subgroups of
 $S_4^+$ are as follows:
 \begin{enumerate}
  \item $S_4^+\simeq SO_3^{-1}$.
  \item The quantum orthogonal group $O_2^{-1}$.
  \item The quantum group $\widehat{D}_{\infty}$, the quantum dual of the infinite dihedral group.
  \item The symmetric group $S_4$ and its subgroups.
  \item The quantum group $S_4^{\tau}$, the unique non-trivial twist
   of $S_4$.
  \item The quantum group $A_5^{\tau}$, the unique non-trivial twist
  of the alternating group $A_5$.
  \item The quantum group $D_n^\tau$, $n$ even and $n\geq 6$, the unique non-trivial
  twist of the dihedral group of order $2n$.
\item The quantum group $DC_n^\tau$ of order $4n$, $n\geq 2$, a pseudo-twist of the dicyclic group of order $4n$.
  \item The quantum group $\widehat{D}_n$, $n\geq 3$, the quantum dual of the dihedral group of order $2n$.
 \end{enumerate}
 \end{theorem}

 The first step in the proof is to show that
 $A_s(4)$ is in fact isomorphic with $C(SO_3^{-1})$,
 the latter being the quotient of $C(SU_3^{-1})$
 by the relations making the fundamental matrix orthogonal.
 Then one shows that $C(SO_3^{-1})$ is a twist of $C(SO_3)$
 (recall that, in contrast, $C(SU_2^{-1})$ is not a twist of
 $C(SU_2)$).
 Then one uses twisting techniques
 to show that the quantum subgroups diagonally containing
 the Klein subgroup of $SO_3^{-1}$ correspond to
 twists of subgroups of $SO_3$ containing
 the diagonal subgroup.
 The remaining cases are examined by using case-by-case arguments. 

The existence and uniqueness of various quantum groups
 in the theorem follow from work of several authors, including
 Davydov, Etingof and Gelaki, Kac and Paljutkin, Masuoka, Nikshych, Vainerman.

 We note that all the quantum groups occurring in the
 theorem were already known, and that the classification
 has lots of similarities with the one for the
 compact subgroups of $SO_3$, which is explained by the twisting
 result.

 A direct consequence of the classification theorem is
 the following result.

 \begin{theorem}
 The classical symmetric group $S_4$ is maximal as a compact quantum subgroup of  the quantum permutation
 group $S_4^+$.
\end{theorem}

We conjecture that for any $n$, the classical symmetric group $S_n$ is maximal in the quantum permutation
 group $S_n^+$.

 The next step is to continue the classification for the next values
 of $n$, say $n=5,6,7$. At this stage we are very far from having complete
 results, or even from having a strategy. We expect that several additional technical difficulties
 will arise. The first one is the non-amenability
 of the discrete quantum group dual to $S_n^+$ if $n \geq 5$, shown in \cite{ba1}.

As a last remark, the results presented here are in connection
with the various $n=4$ results from previous sections. A first
problem here would be to find explicit matrix models for all
quantum groups in the above theorem. Another problem is to
understand the relation with integration results, say via a
systematic study of twisted integration.

\section{Conclusion}

In this paper we have presented several known facts about quantum
permutation groups, most of them being published, or available
from preprints at arxiv.org, and some of them being subject to
papers in preparation. The theory is quite recent: it originates
from Wang's 1998 paper \cite{wa2}, and was basically developed in
the last few years.

The meaning of these investigations might remain quite unclear.
This is indeed the case: the whole subject, with all its possible
interpretations, belongs to the area of mathematical physics, where
everything is by definition quite unclear.

The problem is that the theory is not mature enough for a serious
comparison with results in traditional theoretical physics. Not to
talk about concrete data coming from experimental physics. Wait
and see, and it is most likely that many years will pass before
reaching to the correct technical level.

It is probably instructive here to recall the story of
noncommutative geometry, which is illustrating for the
difficulty of applying mathematical ideas. The theory was
initiated by Connes in the early eighties, with ideas coming from
foliations, groupoids, and the Atiyah-Singer theorem. The high
energy physics motivation was revealed 10 years later, in the
Connes-Lott paper \cite{cl}. This was still quite away from
reasonable numeric results, and the 170 GeV prediction for the
mass of the Higgs boson was obtained 15 more years later, in the
Chamseddine-Connes-Marcolli paper \cite{ccm}.

Back to quantum permutation groups, what we can say for the moment
is that these encode, via a very simple formalism, a few recent results:
\begin{enumerate}
\item Relation with Jones theory. The idea here is to develop a
double approach to the problem, in terms of subfactors and planar
algebras, by using tools from functional analysis and
low-dimensional topology. The first thing to be said is that the
main problems concern the case of non-integer index, and that our
study is quite far away from that (see \cite{bd}). However, in the
simplest case, namely when the index is integer, quantum groups
turn to show up, and we have reasons to believe that quantum
permutation groups have to be investigated first. Among others,
the results presented in this paper are in tune with the recent
work of Bisch and Jones, in the case where the index is integer or
generic. \item Relation with Voiculescu theory. The idea here is
to develop a systematic approach to freeness. Among several
aspects of the theory we have Speicher's approach, in the spirit
of Rota's treatment of classical probability \cite{sp}, and
Biane's representation theory approach \cite{bia}. Our results
make appear a number of key notions of the theory, in the
combinatorial framework. It is our hope that some further random
matrix investigations will lead as well to analytic results.\item
Among other results, we have the pleasing appearance of the Pauli
matrices in connection with the central object of the theory,
namely the algebra $A_s(4)$. Another important fact is the
appearance of Di Francesco's meander determinants in the context
of integration problems. However, regarding both subjects, there
is still a lot of work to be done before understanding why these
objects appear as they do, in the context of quantum permutation
groups.
\end{enumerate}

On the mathematical side, there are a few directions to explore as
well. A number of techniques coming from discrete groups are developed by Vaes, Vergnioux and
their collaborators in \cite{bv}, \cite{bdv}, \cite{va},
\cite{vvv}, \cite{vv}, \cite{ve1}, \cite{ve2}, \cite{ve3}, for various types of discrete quantum groups. Some of these techniques are expected to apply to duals of quantum permutation groups, as a
complement to the above-mentioned considerations.

\end{document}